\documentclass{article}

\usepackage {amsmath, amsthm, amssymb}
\usepackage [top=1.5in, right=1in, bottom=1.5in, left=1in] {geometry}
\usepackage {tikz}
\usepackage {multicol}
\usepackage{hyperref}
\hypersetup{
    colorlinks = true,
    citecolor = blue,
}

\usepackage {libertine}
\usepackage [T1]{fontenc}

\theoremstyle {Theorem}
\newtheorem {thm}  {Theorem}[section]
\newtheorem* {thm*}  {Theorem}
\newtheorem {prop} [thm] {Proposition}
\newtheorem {lem}  [thm]{Lemma}
\newtheorem {cor} [thm] {Corollary}
\newtheorem {ques} [thm] {Question}

\theoremstyle {definition}
\newtheorem {ex} [thm]{Example}
\newtheorem {defn} [thm]{Definition}
\newtheorem {rmk} [thm]{Remark}

\DeclareMathOperator{\ord}{ord}

\title {$q$-Rational and $q$-Real Binomial Coefficients}
\author {John Machacek and Nicholas Ovenhouse}
\date{}

\begin {document}

\maketitle

\begin {abstract}
    We consider $q$-binomial coefficients built from the $q$-rational and $q$-real numbers defined by Morier-Genoud and Ovsienko in terms of continued fractions.
    We establish versions of both the $q$-Pascal identity and the $q$-binomial theorem in this setting.
    These results are then used to find more identities satisfied by the $q$-analogues of Morier-Genoud and Ovsienko, 
    including a Chu--Vandermonde identity and $q$-Gamma function identities.
\end {abstract}

\tableofcontents

\section {Introduction}

The classical ``\emph{$q$-integer}'' is the polynomial $[n]_q = 1 + q + q^2 + \cdots + q^{n-1}$. Since this is also equal to $\frac{1-q^n}{1-q}$,
it is also common to define the $q$-analogue of a non-integer number $\alpha$ by the expression $\frac{1-q^\alpha}{1-q}$. Recently,
Morier-Genoud and Ovsienko have defined a different $q$-analogue of rational~\cite{mgo} and real numbers~\cite{mgo_22} by using $q$-deformations of their
continued fraction expressions. For $\alpha \in \Bbb{R}$, we will use $[\alpha]_q$ to denote the $q$-analog from \cite{mgo,mgo_22}. 
In the case of $\alpha \in \Bbb{Q}$, this definition produces a rational function of $q$, unlike the more common $q$-number $\frac{1-q^\alpha}{1-q}$.
More generally, for $\alpha \in \Bbb{R}$ the $q$-analogue is a Laurent series in $q$ with integer coefficients.

The main goal of this paper is to demonstrate that many definitions involving the $q$-numbers $\frac{1-q^\alpha}{1-q}$
have alternative versions which use the Morier-Genoud--Ovsienko $q$-analogs $[\alpha]_q$, and that many of the results and identities involving $\frac{1-q^\alpha}{1-q}$
are still true for the versions which use $[\alpha]_q$. Our main example will be to generalize the $q$-binomial coefficients $\binom{n}{k}_q$
to the case when $n$ is not an integer, and to see that many of the expected properties and identities hold in this new case.

The binomial coefficients can be written as $\binom{n}{k} = \frac{n(n-1)(n-2) \cdots (n-k+1)}{k!}$, and this of course makes sense even
when $n$ is not an integer. Likewise, we can define $q$-analogues 
\[ \binom{\alpha}{k}_q := \frac{[\alpha]_q [\alpha-1]_q \cdots [\alpha-k+1]_q}{[k]_q!} \] 
which will be Laurent series in $q$ (and in fact a rational function if $\alpha \in \Bbb{Q}$). Again,
we stress that this sort of definition has been made before (see e.g. \cite{kac_book}) 
using the expression $\frac{(1-q^\alpha)(1-q^{\alpha-1}) \cdots (1-q^{\alpha-k+1})}{(1-q)(1-q^2) \cdots (1-q^k)}$,
but the latter is not a rational function nor even a Laurent series in $q$.
This latter version of $q$-binomial coefficient satisfies many nice properties.
They satisfy $q$-Pascal identities and $q$-analogues of the binomial theorem. We will show that our $\binom{\alpha}{k}_q$, using
the $q$-numbers $[\alpha]_q$, satisfy all of the same nice identities.

The main idea which permits the generalizations mentioned above is to replace the expression $q^\alpha$ with some expression which is a 
Laurent series in $q$ (or a rational function for $\alpha \in \Bbb{Q}$), but which retains many of the nice properties which make 
$\frac{1-q^\alpha}{1-q}$ a good $q$-analogue of the number $\alpha$.
Note that when $n$ is an integer, then $[n+1]_q - [n]_q = q^n$. We suggest that for real numbers $\alpha$, the function
$\{\alpha\}_q := [\alpha+1]_q - [\alpha]_q$ is the appropriate substitute for $q^\alpha$. The first observation which indicates that
this is a good idea is the identity $[\alpha]_q = \frac{1-\{\alpha\}_q}{1-q}$. This observation, though very simple, seems to have many
nice and interesting consequences, allowing one to use $\{\alpha\}_q$ in place of $q^\alpha$ in many existing $q$-analogues.
As mentioned above, we will focus in this paper on the $q$-binomial coefficients $\binom{\alpha}{k}_q$, and see that many identities
involving them hold in this new setting after replacing $q^\alpha$ by $\{\alpha\}_q$.

\section {$q$-Rational and $q$-Real Numbers} \label{sec:q_rational_definition}

We will recall the definition of $q$-deformed rational numbers from \cite{mgo}. 

\begin {defn} \label{def:q_rational}
    If $\alpha > 1$ has continued fraction expansion
    \[ 
        \alpha = [a_1,a_2, \dots, a_{2m}] = 
        a_1 + \cfrac{1}{
            a_2 + \cfrac{1}{
                \ddots \, + \, \cfrac{1}{
                    a_{2m}
                } 
            }
        } 
    \]
    then we define the $q$-analogue of $\alpha$ as the rational function $[\alpha]_q$ given by 
    the $q$-deformed continued fraction:
    \[ 
        [\alpha]_q := 
        [a_1]_q + \cfrac{q^{a_1}}{
            [a_2]_{q^{-1}} + \cfrac{q^{-a_2}}{
                [a_3]_q + \cfrac{q^{a_3}}{
                    [a_4]_{q^{-1}} + \cfrac{q^{-a_4}}{
                        \ddots \, + \, \cfrac{q^{a_{2m-1}}}{
                            [a_{2m}]_{q^{-1}}
                        } 
                    }
                }
            }
        } 
    \]
\end {defn}

\bigskip

\begin {ex} \label{ex:q_rational_example}
    The continued fraction for $\alpha = \frac{52}{23}$ is $[2,3,1,5]$, so its $q$-analogue is
    \[ 
        \left[ \frac{52}{23} \right]_q = [2]_q + \cfrac{q^2}{
                                                     [3]_{q^{-1}} + \cfrac{q^{-3}}{
                                                                        1 + \cfrac{q}{
                                                                                [5]_{q^{-1}}
                                                                            }
                                                                    }
                                                 }
        = \frac{1 + 3q + 5q^2 + 7q^3 + 8q^4 + 8q^5 + 7q^6 + 6q^7 + 4q^8 + 2q^9 + q^{10}}{1 + 2q + 3q^2 + 4q^3 + 4q^4 + 3q^5 + 3q^6 + 2q^7 + q^8}
    \]
\end {ex}

In~\cite{mgo_22} the authors obtain a formal series $q$-analog $[\alpha]_q$ for any positive real number $\alpha$ by first taking a sequence of 
rational numbers $(\alpha_j)_{j \geq 0}$ converging to $\alpha$. Next the formal power series $[\alpha]_q$ is defined to be the limit of $([\alpha_j]_q)_{j \geq 0}$ 
(which is shown to be well-defined). 

\begin {thm} [\cite{mgo_22}, Theorem 1] \label{thm:q_real}
    Let $\alpha \in \Bbb{R} \setminus \Bbb{Q}$, and let $\alpha_1,\alpha_2,\alpha_3, \dots$ be a sequence of rationals with $\lim_{n \to \infty} \alpha_n = \alpha$.
    Then the limit $\lim_{n \to \infty} [\alpha_n]_q$ exists as a formal Laurent series with integer coefficients (i.e. each coefficient eventually stabilizes).
\end {thm}

\begin {defn}
    For irrational $\alpha \in \Bbb{R}$, the $q$-analogue $[\alpha]_q \in \Bbb{Z}((q))$ is defined as the Laurent series guaranteed by
    Theorem \ref{thm:q_real}. That is, for any sequence of rationals $\alpha_n \to \alpha$, if $[\alpha_n]_q = \sum_k \varkappa_k^{(n)} q^k$,
    then define
    \[ [\alpha]_q := \sum_k \varkappa_k q^k, \quad \quad \text{ where } \varkappa_k := \lim_{n \to \infty} \varkappa_k^{(n)} \]
\end {defn}

\bigskip

\begin {ex}
    The continued fraction expansion of $\pi/2$ begins with $[1,1,1,3,31,1,145,\dots]$. We can therefore approximate it
    by the sequence of rationals
    \[ \frac{3}{2}, ~ \frac{11}{7}, ~ \frac{344}{219}, ~ \frac{355}{226}, ~ \dots \]
    The series expansions of these $q$-rationals are given by
    \begin {align*}
        \left[ \frac{3}{2} \right]_q     &= 1 + q^2 - q^3 + q^4 - q^5 + q^6 - q^7 + q^8 - q^9 + \cdots \\
        \left[ \frac{11}{7} \right]_q    &= 1 + q^2 - q^3 + q^4 - 2 q^6 + 4 q^7 - 5 q^8 + 4 q^9 - 7q^{11} + \cdots \\
        \left[ \frac{344}{219} \right]_q &= 1 + q^2 - q^3 + q^4 - q^5 + 2 q^6 - 3 q^7 + 3 q^8 - 4 q^9 + 6q^{10} - 7q^{11} + \cdots
    \end {align*}
    The series $\left[ \frac{344}{219} \right]_q$ already agrees with $\left[ \frac{\pi}{2} \right]_q$ to much higher order than what is shown above.
\end {ex}

\bigskip

\begin {prop}[\cite{mgo, mgo_22}] \label{prop:add_1}
    For a real number $\alpha > 1$, we have
    \[ \left[ \alpha + 1 \right]_q = 1 + q \left[ \alpha \right]_q \]
\end {prop}

Proposition \ref{prop:add_1} follows easily from Definition \ref{def:q_rational}, and can be used to define the $q$-analog for $\alpha \leq 1$.
That is, we simply define $[\alpha - 1]_q := \frac{1}{q} \left( [\alpha]_q - 1 \right)$.
We then obtain the following corollary by iteration, which shows that $q$-rational numbers need not always be power series in $q$, but can in fact be Laurent series.

\begin {cor} \label{cor:add_integers}
    For a real number $\alpha \in \Bbb{R}$, and a positive integer $n$, we have
    \[ [\alpha+n]_q = [n]_q + q^n [\alpha]_q \]
    \[ [\alpha-n]_q = \frac{1}{q^n} \Big( [\alpha]_q - [n]_q \Big) \] 
\end {cor}

\bigskip

For integers, the $q$-analogue $[n]_q = 1+q+q^2+\cdots+q^{n-1}$ is equal to $\frac{1-q^n}{1-q}$, and for this reason it 
is common to take the $q$-analogue of an arbitrary real number $\alpha$ to be the expression $\frac{1-q^{\alpha}}{1-q}$.
The $q$-rational numbers from Definition \ref{def:q_rational} are quite different from this expression, but they have many
interesting properties and connections to other areas of mathematics, some of which we will briefly mention now.

The numerators and denominators of $q$-rationals have interesting combinatorial interpretations as certain generating functions.
Let $[a_1,\dots,a_{2m}] = \frac{r}{s}$, with $q$-version given by the rational function $\left[ \frac{r}{s} \right]_q = \frac{\mathcal{R}(q)}{\mathcal{S}(q)}$.
There is a certain poset $F$ (called a \emph{fence}), where the numbers $a_1,\dots,a_{2m}$ determine the number of up and down steps in the Hasse diagram.
Then the polynomial $\mathcal{R}(q)$ is the rank generating function of the lattice $J(F)$, of order ideals of $F$. That is,
\[ \mathcal{R}(q) = \sum_{I \in J(F)} q^{|I|} \]
A version of this statement appears in \cite{mgo} (Theorem 4). These polynomials were also studied from this fence poset perspective in \cite{mss}.
A different, yet equivalent, perspective is the following. There is a \emph{border strip} skew Young diagram $\lambda / \mu$, whose shape is determined
by the numbers $a_1,\dots,a_{2m}$. It is obtained by rotating the Hasse diagram of the fence poset $F$ by $45^\circ$ counter-clockwise, and then replacing
each vertex with a box. Then there is a simple bijection between order ideals of $F$ and north-east lattice paths on $\lambda / \mu$, where the order ideal
is the set of boxes underneath the lattice path. If we let $L(\lambda / \mu)$ be the set of north-east lattice paths on the shape $\lambda / \mu$, and
for a path $p$, let $|p|$ be the number of boxes under the path, the statement above can be re-phrased as
\[ \mathcal{R}(q) = \sum_{p \in L(\lambda / \mu)} q^{|p|} \]
This interpretation, in terms of lattice paths and skew Young diagrams, was used in \cite{ovenhouse_21} to give another combinatorial
meaning of the polynomials $\mathcal{R}(q)$ as counting the sizes of certain varieties over the finite field with $q$ elements.
In particular, $q^{|\mu|} \mathcal{R}(q)$ is the number of $\Bbb{F}_q$-points in a union of Schubert cells in some Grassmannian,
where the union is over the Schubert cells indexed by partitions $\nu$ with $\mu \leq \nu \leq \lambda$.

The $q$-rationals also have a significant connection to cluster algebras. Cluster algebras of ``type A'' can be realized as the
homogeneous coordinate rings of Grassmannians $\mathrm{Gr}_2(n)$, and also as the ring of functions on Penner's \emph{decorated Teichm\"{u}ller space}
of an ideal polygon in the hyperbolic plane. In both cases, the cluster combinatorics are determined by triangulations of a polygon. 
Any triangulation of a polygon, and a choice of a diagonal not in that triangulation, determines a continued fraction $[a_1,\dots,a_{2m}]$.
It was explained in \cite{mgo} (appendix B.2) that $\mathcal{R}(q)$ is a specialization of the $F$-polynomial of the corresponding cluster variable.

A nice exposition of some of the combinatorial formulas mentioned above, and the connection with cluster algebras, is given in \cite{claussen}.
The $q$-rationals are also related to Jones polynomials of certain knots. This is explained in \cite{mgo} (appendix A). This is related to 
work of Lee and Schiffler \cite{ls_19}, who showed the Jones polynomials of certain knots have a cluster algebra interpretation.

Some instances where $q$-real are used include the deformation of modular groups in \cite{mgl} and also of Conway--Coxeter friezes in \cite{conwayCox}.
The $q$-real numbers are (Laurent) series, and their convergence properties have been studied in \cite{radius}. Also in \cite{mgl}, it was shown
that the $q$-deformations of quadratic irrationals have a particularly nice closed-form algebraic expression.

\bigskip

\section {A Substitute for $q^\alpha$}

It is our hope that we might convince the reader not only that $q$-rationals and $q$-reals are interesting for the reasons mentioned in the previous section, but
also that they satisfy many nice identities analogous to those satisfied by the more commonly used expressions $\frac{1-q^{\alpha}}{1-q}$. 
Many identities involving these $q$-analogues will inevitably involve the expression $q^\alpha$, which is certainly not a rational function.
In this section, we define a Laurent series denoted $\{\alpha\}_q$ (which is a rational function for $\alpha \in \Bbb{Q}$) which we suggest is a good replacement for $q^\alpha$. 
For the remainder of the paper, we will demonstrate many examples of identities satisfied by the $q$-rational and $q$-real
numbers which resemble well-known identities, but with $q^\alpha$ replaced by $\{\alpha\}_q$.

Note the simple fact that when $n$ is an integer, $q^n = [n+1]_q - [n]_q$. We propose that this forward difference is the ``correct'' substitute
for $q^{\alpha}$ when $\alpha$ is not an integer.

\bigskip

\begin {defn}
    For $\alpha \in \Bbb{R}$, define $\{\alpha\}_q := [\alpha+1]_q - [\alpha]_q$.
\end {defn}

\bigskip

\begin {ex} \label{ex:q^alpha_examples}
    Here are some examples of $\{\alpha\}_q$ for rational values of $\alpha$: 
    \[ 
        \left\{ \frac{1}{2} \right\}_q = \frac{1+q^2}{1+q}, \quad 
        \left\{ \frac{5}{3} \right\}_q = q \frac{1+q^2+q^3}{1+q+q^2}, \quad 
        \left\{ \frac{25}{7} \right\}_q = q^3 \frac{1 + q + 2q^2 + q^3 + q^4 + q^5}{1 + 2q + 2q^2 + q^3 + q^4} 
    \]
    Note that in all these examples, there is a factor of $q^{\lfloor \alpha \rfloor}$. This is always the case, and this follows from 
    Theorem \ref{thm:order_of_alpha} and Proposition \ref{prop:q^alpha_properties}(c) below.
\end {ex}

\bigskip

\begin {ex}
    For any positive integer $n$, we have
    \[ 
        \left\{ \frac{1}{n} \right\}_q = \frac{1 + q + q^2 + \cdots + q^{n-2} + q^n}{1 + q + q^2 + \cdots + q^{n-1}} 
        = 1 + (q-1) \frac{q^{n-1}}{[n]_q}
    \]
    Expanded as a Laurent series, we have $\left\{\frac{1}{n}\right\}_q = 1 - q^{n-1} + 2q^n + \cdots$. 
    In particular, we have $\lim_{n \to \infty} \left\{\frac{1}{n}\right\}_q = 1$. This limit is to be understood in the sense
    of formal power series, meaning each coefficient eventually stabilizes. We will always work in the ring of formal power series
    $\Bbb{Q}[[q]]$ or Laurent series $\Bbb{Q}((q))$, and limits will always be understood in this sense.
\end {ex}

\bigskip

\begin {ex}
    For an irrational example, we consider $\alpha = \sqrt{2} + 1 \approx 2.41421$. The $q$-analogue $[\sqrt{2}+1]_q$ was computed in \cite{mgo_22} to be
    \[ \left[ \sqrt{2} + 1 \right]_q = 1 + q + q^4 - 2q^6 + q^7 + 4q^8 - 5q^9 - 7q^{10} + 18q^{11} + 7q^{12} - 55q^{13} + \cdots \]
    Using Proposition \ref{prop:q^alpha_properties}(a) below, we can use this expression to compute $\{\sqrt{2}+1\}_q$:
    \[ \left\{ \sqrt{2} + 1 \right\}_q = q^2 - q^4 + q^5 + 2q^6 - 3q^7 - 3q^8 + 9q^9 + 2q^{10} - 25 q^{11} + 11q^{12} + 62q^{13} + \cdots \]
    Note how, just as in Example \ref{ex:q^alpha_examples}, this can be written as $q^{\lfloor \alpha \rfloor}$ times a power series
    with constant term 1. 
\end {ex}

\bigskip

The following summarizes some properties of $\{\alpha\}_q$. We invite the reader to keep in mind the analogy of $\{\alpha\}_q$ with $q^\alpha$,
and to note that parts (c) and (d) correspond to basic rules of exponents $q^{\alpha+n} = q^\alpha q^n$ and $q^{-\alpha} = (q^{-1})^\alpha$.
Furthermore, part (b) shows that the $q$-rationals $[\alpha]_q$ also have a form which resembles $\frac{1-q^\alpha}{1-q}$.

\begin {prop} \label{prop:q^alpha_properties}
    For $\alpha \in \Bbb{R}$ and a positive integer $n$,
    \begin {enumerate}
        \begin {multicols}{2}
        \item[$(a)$] $\{\alpha\}_q = 1 + (q-1)[\alpha]_q$
        \item[$(b)$] $\displaystyle \frac{1 - \{\alpha\}_q}{1-q} = [\alpha]_q$
        \item[$(c)$] $\{\alpha+n\}_q = q^n \{\alpha\}_q$
        \columnbreak
        \item[$(d)$] If $\alpha \in \Bbb{Q}$, then $\{-\alpha\}_q = \{\alpha\}_{q^{-1}}$ 
        \item[$(e)$] $\displaystyle \{\alpha\}_q = \frac{[\alpha+n]_q - [\alpha]_q}{[n]_q}$
        \end {multicols}
    \end {enumerate}
\end {prop}
\begin {proof}
    Part $(a)$ follows immediately from Proposition \ref{prop:add_1}. Part $(b)$ follows from $(a)$ by algebraic manipulation.

    For part $(c)$, note that by a combination of part $(a)$ and Proposition \ref{prop:add_1}, we have 
    \begin {align*} 
        \{\alpha+1\}_q &= 1 + (q-1)[\alpha+1]_q                \tag{part (a)}              \\
                       &= 1 + (q-1)(1 + q[\alpha]_q)           \tag{Prop \ref{prop:add_1}} \\
                       &= q \left( 1 + (q-1)[\alpha]_q \right)                             \\
                       &= q \{\alpha\}_q                       \tag{part (a)}
    \end {align*}
    The general case for $\{\alpha+n\}_q$ then follows by induction.

    For part $(d)$, we use the fact  that $[-\alpha]_q = -q^{-1}[\alpha]_{q^{-1}}$ (Proposition 2.8 from \cite{mgl}).
    Then we have
    \begin {align*}
        \{-\alpha\}_q &= [-(\alpha-1)]_q - [-\alpha]_q                                                                          \\
                      &= q^{-1} \Big( [\alpha]_{q^{-1}} - [\alpha-1]_{q^{-1}} \Big)      \tag{\cite{mgl}, Prop 2.8}             \\
                      &= q^{-1} \Big( [\alpha]_{q^{-1}} - q([\alpha]_{q^{-1}} - 1) \Big) \tag{Corollary \ref{cor:add_integers}} \\
                      &= 1 + (q^{-1} - 1)[\alpha]_{q^{-1}}                                                                      \\
                      &= \{\alpha\}_{q^{-1}}                                             \tag{part (a)}
    \end {align*}

    Part $(e)$ is a simple calculation using Corollary \ref{cor:add_integers} and the fact that $[n]_q = \frac{1-q^n}{1-q}$.
\end {proof}

\bigskip

\begin {rmk} \label{rmk:alpha^k}
    By Corollary \ref{cor:add_integers}, we have $[\alpha + k]_q = [k]_q + q^k [\alpha]_q$, and part $(e)$ of Proposition \ref{prop:q^alpha_properties} is
    saying we also have $[\alpha+k]_q = [\alpha]_q + \{\alpha\}_q [k]_q$, so there is some symmetry in the roles of $\alpha$ and $k$. 
    However, $[\alpha+\beta]_q \neq [\alpha]_q + \{\alpha\}_q [\beta]_q$ in general.
    It is only true if either $\alpha$ or $\beta$ is an integer. For example, when $\alpha=\beta=1/2$, we have $[1]_q = 1$, 
    but on the other hand, $\left[ \frac{1}{2} \right]_q + \left\{ \frac{1}{2} \right\}_q \left[ \frac{1}{2} \right]_q = q \frac{2 + q + q^2}{1 + 2q + q^2}$.
\end {rmk}

\bigskip

\begin {rmk}
    In part $(c)$ of Proposition \ref{prop:q^alpha_properties}, it is important that $n$ is an integer. The more general version of this statement is \emph{not} true.
    That is, $\{\alpha+\beta\}_q \neq \{\alpha\}_q \{\beta\}_q$ in general. For example, when $\alpha = \frac{1}{2}$ and $\beta = \frac{3}{2}$,
    we have $\{2\}_q = q^2$, but $\left\{ \frac{1}{2} \right\}_q \left\{ \frac{3}{2} \right\}_q = q \frac{1+2q^2+q^4}{1+2q+q^2}$.
\end {rmk}

\bigskip

\begin {rmk}
    As mentioned in \cite{mgl}, part (d) of Proposition \ref{prop:q^alpha_properties} only makes sense for $\alpha \in \Bbb{Q}$.
\end {rmk}

\bigskip

\section {$q$-Rational and $q$-Real Binomial Coefficients}

\bigskip

The $q$-integers $[n]_q = 1+q+ \cdots + q^{n-1}$ are the building blocks of many other $q$-analogues. Perhaps the most well-known are
the $q$-factorial $[n]_q ! = [n]_q [n-1]_q \cdots [2]_q [1]_q$ and the $q$-binomial coefficients $\binom{n}{k}_q = \frac{[n]_q!}{[k]_q![n-k]_q!}$.
Another very useful notation in the theory of $q$-analogues is the $q$-Pochhammer symbol
\[ (x;q)_n := (1-x)(1-qx)(1-q^2x) \cdots (1-q^{n-1}x) \]
\[ (x;q)_\infty := \lim_{n \to \infty} (x;q)_n = \prod_{k=0}^\infty (1-q^kx) \]
The $q$-factorial and $q$-binomial coefficients
can both be expressed in terms of the $q$-Pochhammer symbol by
\[ [n]_q! = \frac{(q;q)_n}{(1-q)^n} \quad \text{ and } \quad \binom{n}{k}_q = \frac{(q;q)_n}{(q;q)_k (q;q)_{n-k}} \]
We generalize the $q$-binomial coefficient to the case when $n = \alpha \in \Bbb{R}$ is not an integer by the following formula.

\begin {defn}
    For $\alpha \in \Bbb{R}$, define the $q$-binomial coefficient $\binom{\alpha}{k}_q$ as
    \[ \binom{\alpha}{k}_q := \frac{[\alpha]_q [\alpha-1]_q \cdots [\alpha-k+1]_q}{[k]_q!} \]
    where the factors in the numerator are $q$-rational or $q$-real numbers. Note that $\binom{\alpha}{k}_q$ is
    a Laurent series (with integer coefficients), which happens to be a rational function when $\alpha \in \Bbb{Q}$.
\end {defn}

\bigskip

\begin {ex}
    $\displaystyle \binom{5/3}{3}_q = \frac{[5/3]_q [2/3]_q [-1/3]_q}{[3]_q [2]_q [1]_q} = \frac{-(1 + q + 2q^2 + q^3)}{1 + 4q + 10q^2 + 16q^3 + 19q^4 + 16q^5 + 10q^6 + 4q^7 + q^8}$
\end {ex}

\bigskip

These generalized $q$-binomial coefficients satisfy the usual $q$-Pascal identities, using $\{\alpha\}_q$ in place of $q^\alpha$ where appropriate.

\begin {prop} \label{prop:q_Pascal}
    For a real number $\alpha \in \Bbb{R}$ and integer $k \geq 0$, we have
    \begin {enumerate}
        \item[$(a)$] $\displaystyle \binom{\alpha}{k}_q = q^k \binom{\alpha-1}{k} + \binom{\alpha-1}{k-1}$
        \item[$(b)$] $\displaystyle \binom{\alpha}{k}_q = \binom{\alpha-1}{k} + \{\alpha-k\}_q \binom{\alpha-1}{k-1}$
    \end {enumerate}
\end {prop}
\begin {proof}
    Thanks to Corollary \ref{cor:add_integers}, Proposition \ref{prop:q^alpha_properties}(e), and Remark \ref{rmk:alpha^k}, 
    the proof is essentially the same as for the classical $q$-binomial coefficients (i.e. the case when $\alpha$ is an integer).

    $(a)$ Look at the right-hand side:
    \begin {align*}
        q^k \binom{\alpha-1}{k} + \binom{\alpha-1}{k-1} 
        &= \frac{q^k[\alpha-1]_q [\alpha-2]_q \cdots [\alpha-k]_q}{[k]_q!} + \frac{[\alpha-1]_q[\alpha-2]_q \cdots [\alpha-k+1]_q}{[k-1]_q!} \\
        &= \frac{[\alpha-1]_q[\alpha-2]_q \cdots [\alpha-k+1]_q}{[k]_q!} \Big( q^k [\alpha-k]_q + [k]_q \Big)
    \end {align*}
    By Corollary \ref{cor:add_integers}, the expression in parentheses is equal to $[\alpha]_q$, and this gives the result.

    \medskip

    $(b)$ Again, look at the right-hand side and combine:
    \begin {align*}
        \binom{\alpha-1}{k} + \{\alpha-k\}_q \binom{\alpha-1}{k-1} 
        &= \frac{[\alpha-1]_q [\alpha-2]_q \cdots [\alpha-k]_q}{[k]_q!} + \frac{[\alpha-1]_q[\alpha-2]_q \cdots [\alpha-k+1]_q \{\alpha-k\}_q}{[k-1]_q!} \\
        &= \frac{[\alpha-1]_q[\alpha-2]_q \cdots [\alpha-k+1]_q}{[k]_q!} \Big( [\alpha-k]_q + \{\alpha-k\}_q [k]_q \Big)
    \end {align*}
    By Remark \ref{rmk:alpha^k}, the expression in parentheses is equal to $[\alpha]_q$.
\end {proof}

\bigskip

Now we give several other formulas for $\binom{\alpha}{k}_q$, some of which will be useful later.

\begin {prop} \label{prop:alternate_binom_formula}
    The $q$-binomial coefficients have the following alternate formulations: for $\alpha \in \Bbb{R}$,
    \begin {enumerate}
        \begin {multicols}{2}
        \item[$(a)$] $\displaystyle \binom{\alpha}{k}_q = \prod_{i=0}^{k-1} \frac{[\alpha]_q - [i]_q}{[k]_q - [i]_q}$ 
        \item[$(b)$] $\displaystyle \binom{\alpha}{k}_q = \frac{1}{q^{\binom{k}{2}} [k]_q!} \prod_{i=0}^{k-1} \left( [\alpha]_q - [i]_q \right)$
        \item[$(c)$] $\displaystyle \binom{\alpha}{k}_q = \frac{(-\{\alpha\}_q)^k}{q^{\binom{k}{2}}} \cdot \frac{\left( \{\alpha\}_q^{-1}; q\right)_k}{(q;q)_k}$
        \columnbreak
        \item[$(d)$] $\displaystyle \binom{\alpha}{k}_q = \frac{(\{\alpha\}_q; q^{-1})_k}{(q;q)_k}$
        \item[$(e)$] $\displaystyle \binom{\alpha+k-1}{k}_q = \frac{(\{\alpha\}_q; q)_k}{(q;q)_k}$ 
        \end {multicols}
    \end {enumerate}
\end {prop}
\begin {proof}
    Parts $(a)$ and $(b)$ follow directly from Corollary \ref{cor:add_integers}. 

    To verify $(d)$, use Proposition \ref{prop:q^alpha_properties}, parts (b) and (c), to write
    \begin {align*}
        \binom{\alpha}{k}_q &= \frac{(1-\{\alpha\}_q)(1-\{\alpha-1\}_q) \cdots (1 - \{\alpha-k+1\}_q)}{(q;q)_k} \tag{\ref{prop:q^alpha_properties}(b)} \\
                            &= \frac{(1-\{\alpha\}_q)(1-q^{-1}\{\alpha\}_q) \cdots (1 - q^{-(k-1)}\{\alpha\}_q)}{(q;q)_k} \tag{\ref{prop:q^alpha_properties}(c)}
    \end {align*}
    The numerator is equal to $(\{\alpha\}_q; q^{-1})_k$, and part $(d)$ follows. Part $(c)$ is easily obtained from $(d)$ by some algebraic manipulations.
    Part $(e)$ follows from a similar calculation as part $(d)$, using Proposition \ref{prop:q^alpha_properties}(c).
\end {proof}

\bigskip

\begin {rmk}
    Parts $(c)$, $(d)$, and $(e)$ of Proposition \ref{prop:alternate_binom_formula} support the claim that $\{\alpha\}_q$ is a good substitute for $q^\alpha$, 
    as they resemble the corresponding formulas when $\alpha = n$ is an integer:
    \[ \binom{n}{k}_q = (-1)^k q^{kn-\binom{k}{2}} \frac{(q^{-n};q)_k}{(q;q)_k} = \frac{(q^n;q^{-1})_k}{(q;q)_k} = \frac{(q^{n-k};q)_k}{(q;q)_k} \]
\end {rmk}

\bigskip

For a (formal) Laurent series $f(q) = \sum_{k} c_k q^k$ we let $\ord(f)$ denote the \emph{order} of $f$ which is the minimal value of $k$ such that $c_k \neq 0$.
The order is a valuation on the ring of Laurent series.
We now give a rephrasing of a theorem from~\cite{mgo_22}.

\begin {thm} [\cite{mgo_22}, Theorem 2] \label{thm:order_of_alpha}
    For $\alpha \in \Bbb{R}$, if $0 < \alpha < 1$ we have $\ord([\alpha]_q) \geq 1$ and otherwise
    \[ 
        \ord([\alpha]_q) = \begin{cases} 
            0      & \text{ if } ~ \alpha \geq 1  \\ 
            \infty & \text{ if } ~ \alpha = 0     \\ 
            N      & \text{ if } ~ N \leq \alpha < N + 1 \text{ and } N \in \Bbb{Z}_{<0}
        \end{cases}
    \]
gives the order of $[\alpha]_q$ exactly.
\label{thm:order}
\end {thm}

We note that when $\alpha < 0$, we have simply that $\ord([\alpha]_q) = \lfloor \alpha \rfloor$.
In the case that $\alpha \geq 0$, we see that the order falls into one of only three possible cases.
In the case that $0 < \alpha < 1$ we only have a bound of the order, but this bound will be sufficient for our purposes.
It will be useful to know the order of the Laurent series given by $q$-binomial coefficients to be able to ensure certain combinations and evaluations are well-defined.
To reduce the number of cases the following lemma does not condsider ordinary integer $q$-binomial coefficeints $\binom{n}{k}$ for a nonegative integer $n$ 
since these are always either a polynomial with nonzero constant term or idenitically zero (this latter case happens when $k > n$).

\begin{lem} \label{lem:order}
    Consider $\alpha \in \Bbb{R} \setminus \Bbb{Z}_{\geq 0}$ and $k \in \Bbb{Z}_{\geq 0}$.
Let $\beta = \alpha - \lfloor \alpha \rfloor$ and set $b = \ord([\beta]_q)$.
    Assume that $N \in \Bbb{Z}$ with $N \leq \alpha < N+1$ so that $\lfloor \alpha \rfloor = N$, then
    \[ \ord\left(\binom{\alpha}{k}_q \right)  = \begin{cases} 0 & k \leq N\\ 
b - \binom{k-N}{2} & k > N \text{ and } N  \geq 0 \\
Nk - \binom{k}{2} & N < 0. \\ \end{cases}\]
\end{lem}
\begin{proof}
    Since $[k]_q!$ is a polynomial with constant term equal to $1$ it follows that $\displaystyle \frac{1}{[k]_q!}$ is a power series with constant term equal to $1$.
    Thus the order of $\binom{\alpha}{k}_q$ is completely determined by the order of the falling factorial $[\alpha]_q [\alpha-1]_q \cdots [\alpha-k+1]_q$ in the numerator.
    First consider the case that $N < 0$, then $\ord([\alpha-j]_q) = N - j$ and
    \[\ord([\alpha]_q [\alpha-1]_q \dots [\alpha - k +1]_q) = \sum_{j=0}^{k-1} (N - j) = Nk - \binom{k}{2}.\]

If $k \leq N$, then $N \geq 0$ and the order will be $0$ as each term in the falling factorial will have order equal to $0$ by Theorem~\ref{thm:order}.
   Last consider the case $N \geq 0$ with $k > N$, then
    \begin{align*}
    \ord([\alpha]_q [\alpha-1]_q \dots [\alpha - k +1]_q) &= \ord([\alpha-N]_q) + \ord([\alpha-N-1]_q \dots [\alpha - k +1]_q) \\
&= \ord([\beta]_q) + \ord([\alpha-N-1]_q \dots [\alpha - k +1]_q) \\
    &= b - \sum_{j=1}^{k-N-1} j\\
 &= b - \binom{k-N}{2}
    \end{align*}
    and the lemma is proven.
\end{proof}

\bigskip

\section {Combinatorial Intepretation}

In Section \ref{sec:q_rational_definition}, we mentioned some combinatorial interpretations of $q$-rational numbers appearing in the literature.
Using these, it is easy to give some combinatorial interpretations of the $q$-rational binomial coefficients.
In this section, we will describe this combinatorial model, and give some examples.

Let $\alpha \in \Bbb{Q}$ with continued fraction $\alpha = [a_1,a_2,\dots,a_{2m}]$. Define a planar graph $G_\alpha$
built out of squares, where each square is attached to the previous on the right or above. Reading the sequence of
``right''s and ``up''s should give $U^{a_1-1} R^{a_2} U^{a_3} R^{a_4} \cdots U^{a_{2m-1}} R^{a_{2m}-1}$. In other words,
the number of consecutive ``up'' and ``right'' steps are given by the continued fraction coefficients (with the first and last being off-by-one).
These graphs $G_\alpha$ are often called ``\emph{snake graphs}'' in the literature (e.g. in \cite{propp} and \cite{cs_18}).

\bigskip

\begin {ex}
    Recall Example \ref{ex:q_rational_example}, where we considered $\frac{52}{23} = [2,3,1,5]$. The corresponding graph has ``up'' and ``right'' steps
    in the sequence $UR^3UR^4$. The graph $G_{52/23}$ is shown in figure \ref{fig:snake_example}.

    \begin {figure}
    \centering
    \begin {tikzpicture}[scale=0.7]
        \draw (0,0) -- (1,0) -- (1,1) -- (4,1) -- (4,2) -- (8,2) -- (8,3) -- (3,3) -- (3,2) -- (0,2) -- cycle;
        \draw (0,1) -- (1,1) -- (1,2);
        \draw (2,1) -- (2,2);
        \draw (3,1) -- (3,2) -- (4,2) -- (4,3);
        \draw (5,2) -- (5,3);
        \draw (6,2) -- (6,3);
        \draw (7,2) -- (7,3);
    \end {tikzpicture}
    \caption {The snake graph $G_\alpha$ for $\alpha = \frac{52}{23}$.}
    \label {fig:snake_example}
    \end {figure}
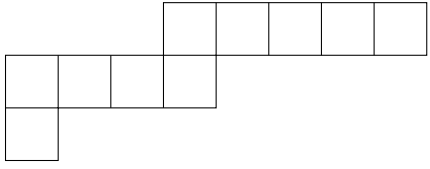
\end {ex}

\bigskip

As mentioned in Section \ref{sec:q_rational_definition}, if $[\alpha]_q = \frac{\mathcal{R}(q)}{\mathcal{S}(q)}$ (with $\alpha > 1$), then the numerator $\mathcal{R}(q)$
is the generating function for an area statistic on north-east lattice paths in $G_\alpha$. More specifically, if $\mathcal{L}(G_\alpha)$
is the set of north-east lattice paths in $G_\alpha$ from the bottom-left to the top-right corner, with $|p|$ denoting the area (i.e. the number of boxes)
underneath the path $p$, then $\mathcal{R}(q) = \sum_{p \in \mathcal{L}(G_\alpha)} q^{|p|}$. This is an equivalent re-statement of Theorem 4 from \cite{mgo},
which was stated in the present form (in terms of the snake graph language) in \cite{claussen} and \cite{ovenhouse_21}. 
Furthermore, the denominator $\mathcal{S}(q)$ has the same interpretation, but for the smaller snake graph obtained from $G_\alpha$ by removing
the initial vertical column of boxes.

\begin {defn} \label{def:L^k}
	For a snake graph $G_\alpha$, let $\mathcal{L}^{(k)}(G_\alpha)$ be the set of $k$-tuples $(p_1,\dots,p_k)$ of north-east lattice paths in $G_\alpha$,
	where each $p_i$  is required to begin with $i-1$ consecutive ``up'' steps. For $p = (p_1,\dots,p_k) \in \mathcal{L}^{(k)}(G_\alpha)$, 
	let $|p| = \sum_{i=1}^k |p_i|$ denote the sum of the areas under the paths.
\end {defn}

\bigskip

\begin {thm}
Suppose $\alpha \in \Bbb{Q}$ with $\alpha > 1$, and let $0 \leq k < \alpha$. The numerator of $\binom{\alpha}{k}_q$ is given by
\[ q^{-\binom{k}{2}} \sum_{p \in \mathcal{L}^{(k)}(G_\alpha)} q^{|p|} \]
\end {thm}
\begin {proof}
Let $[\alpha]_q = \frac{\mathcal{R}(q)}{\mathcal{S}(q)}$ as above.
Note that for $i < \alpha$, $\alpha - i$ will have continued fraction
$[a_1-i, a_2, \dots, a_{2m}]$. Therefore the truncated snake graphs obtained by removing the first column of boxes will be the same for all $\alpha-i$ (for $0 \leq i < k$),
and hence all the $[\alpha-i]_q$ will have the same denominator $\mathcal{S}(q)$. Let $\mathcal{R}_i(q)$ be their numerators, so that $[\alpha-i]_q = \frac{\mathcal{R}_i(q)}{\mathcal{S}(q)}$. Then the $q$-binomial coefficient will be given by
\[ \binom{\alpha}{k}_q = \frac{\mathcal{R}(q) \mathcal{R}_1(q) \mathcal{R}_2(q) \cdots \mathcal{R}_{k-1}(q)}{\mathcal{S}(q)^k [k]_q!} \]
When $i < \alpha$, the graph $G_{\alpha-i}$ is obtained from $G$ by simply removing the first $i$ squares. There is a bijection between $\mathcal{L}(G_{\alpha-i})$ and the set of lattice paths in $G_\alpha$ which begin with $i-1$ up steps (by identifying $G_{\alpha-i}$ with a subgraph of $G_\alpha$).
This bijection is not weight-preserving, but rather the weights differ by a factor $q^{i-1}$ since the corresponding path in $G_\alpha$ goes over $i-1$ extra
boxes compared to its counterpart in $G_{\alpha-i}$. Adding these contributions for all $i=1,2,\dots,k$ gives the extra factor of $q^{-\binom{k}{2}}$.
\end {proof}

\bigskip

\begin {ex}
    The continued fraction for $\frac{5}{2}$ is $[2,2]$. The graph $G_\alpha$ therefore has 3 squares, with up-right sequence $UR$. 
    Let $\mathcal{L}^{(2)}(G_{5/2})$ be the set of ordered pairs of paths $(p_1,p_2)$ in $G_{5/2}$, where $p_2$ starts with an up step.
    The $q$-binomial coefficient $\binom{5/2}{2}_q$ is given by $\frac{1 + 3q + 4q^2 + 4q^3 + 2q^4 + q^5}{(1+q)^3}$. The numerator is
    the sum over $\mathcal{L}^{(2)}(G_{5/2})$ of $q^{|p_1|+|p_2|-1}$. The elements of $\mathcal{L}^{(2)}(G_{5/2})$ are shown in Figure \ref{fig:lattice_path_example}, along with the
    corresponding monomials. In the figure, $p_1$ is in blue, and $p_2$ is in red.

    \begin {figure}
    \centering
    \begin {tikzpicture}[scale=0.5]
        \foreach \x in {0, 4, 8, 12, 16} {
            \foreach \y in {0, 4, 8} {
                \draw (\x,\y) -- (\x+1,\y) -- (\x+1,\y+1) -- (\x+2,\y+1) -- (\x+2,\y+2) -- (\x,\y+2) -- cycle;
                \draw (\x,\y+1) -- (\x+1,\y+1) -- (\x+1,\y+2);
            }
        }

        \foreach \y in {0,4,8} {
            \draw[blue, line width = 2] (0,\y) -- (1,\y) -- (1,\y+1) -- (2,\y+1) -- (2,\y+2);
            \draw[blue, line width = 2] (4,\y) -- (4,\y+1) -- (4+1,\y+1) -- (4+2,\y+1) -- (4+2,\y+2);
            \draw[blue, line width = 2] (8,\y) -- (8+1,\y) -- (8+1,\y+1) -- (8+1,\y+2) -- (8+2,\y+2);
            \draw[blue, line width = 2] (12,\y) -- (12,\y+1) -- (12+1,\y+1) -- (12+1,\y+2) -- (12+2,\y+2);
            \draw[blue, line width = 2] (16,\y) -- (16,\y+2) -- (16+2,\y+2);
        }

        \foreach \x in {0, 4, 8, 12, 16} {
            \draw[red, line width = 2] (\x-0.1,0) -- (\x-0.1,0.9) -- (\x+2.1,0.9) -- (\x+2.1,2);
            \draw[red, line width = 2] (\x-0.1,4) -- (\x-0.1,4+0.9) -- (\x+1.1,4+0.9) -- (\x+1.1,4+2.1) -- (\x+2.1,4+2.1);
            \draw[red, line width = 2] (\x-0.1,8) -- (\x-0.1,8+2.1) -- (\x+2,8+2.1);
        }

        \draw (1,  -0.5) node {$1$};
        \draw (5,  -0.5) node {$q$};
        \draw (9,  -0.5) node {$q$};
        \draw (13, -0.5) node {$q^2$};
        \draw (17, -0.5) node {$q^3$};

        \draw (1,  3.5) node {$q$};
        \draw (5,  3.5) node {$q^2$};
        \draw (9,  3.5) node {$q^2$};
        \draw (13, 3.5) node {$q^3$};
        \draw (17, 3.5) node {$q^4$};

        \draw (1,  7.5) node {$q^2$};
        \draw (5,  7.5) node {$q^3$};
        \draw (9,  7.5) node {$q^3$};
        \draw (13, 7.5) node {$q^4$};
        \draw (17, 7.5) node {$q^5$};
    \end {tikzpicture}
    \caption {An illustration of $\mathcal{L}^{(k)}(G_\alpha)$ for $\alpha = \frac{5}{2}$ and $k=2$. 
                  Below each picture is the monomial $q^{|p|}$, as in Definition \ref{def:L^k}.}
    \label {fig:lattice_path_example}
    \end {figure}
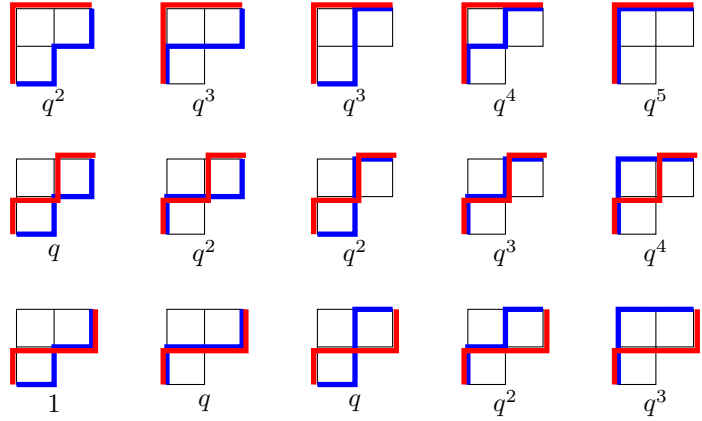
\end {ex}

\bigskip

\section {The $q$-Binomial Theorem}

\bigskip

Recall Newton's generalized binomial theorem: for any $\alpha \in \Bbb{R}$ 
\[ 
    (1+x)^\alpha = \sum_{k=0}^\infty \binom{\alpha}{k} x^k 
    \quad \text{and} \quad 
    \frac{1}{(1-x)^\alpha} = \sum_{k=0}^\infty \binom{\alpha + k - 1}{k} x^k 
\]
When $n$ is a positive integer, there are well-known $q$-versions:
\[ 
    \prod_{k=0}^{n-1} (1 + q^kx) = \sum_{k=0}^n q^{\binom{k}{2}} \binom{n}{k}_q x^k 
    \quad \text{and} \quad
    \prod_{k=0}^{n-1} \frac{1}{1 - q^kx} = \sum_{k=0}^\infty \binom{n + k - 1}{k}_q x^k
\]
These can both be obtained as special cases of the general ``\emph{$q$-binomial theorem}'' (see e.g.~\cite[(1.3.2)]{gasper_book}):
\[ \frac{(ax;q)_\infty}{(x;q)_\infty} = \sum_{k=0}^\infty \frac{(a;q)_k}{(q;q)_k} x^k \]
In particular, substituting $a=q^n$ gives one version, and substituting $x \mapsto -q^nx$ and $a=q^{-n}$ gives the other.
Keeping with our basic philosophy, we will see that using $\{\alpha\}_q$ in place of $q^n$ gives the appropriate
$q$-analogue using $\binom{\alpha}{k}_q$.

\begin {thm} \label{thm:product_formula}
    For $\alpha \in \Bbb{R}$, we have
    \begin {enumerate}
        \item[$(a)$] $\displaystyle \frac{(-x;q)_\infty}{(-\{\alpha\}_q x;q)_\infty} = \prod_{k=0}^\infty \frac{1 + q^kx}{1 + \{\alpha + k\}_q x} = \sum_{k=0}^\infty q^{\binom{k}{2}} \binom{\alpha}{k}_q x^k$
        \item[$(b)$] $\displaystyle \frac{(\{\alpha\}_q x;q)_\infty}{(x;q)_\infty} = \prod_{k=0}^\infty \frac{1 - \{\alpha + k\}_q x}{1 - q^k x} = \sum_{k=0}^\infty \binom{\alpha+k-1}{k}_q x^k$
    \end {enumerate}
\end {thm}
\begin {proof}
    Substitute $a = \{\alpha\}_q$ into the general $q$-binomial theorem to get
    \[ \frac{(\{\alpha\}_qx;q)_\infty}{(x;q)_\infty} = \sum_{k=0}^\infty \frac{(\{\alpha\}_q;q)_k}{(q;q)_k} x^k  \]
    By Proposition \ref{prop:alternate_binom_formula}(e), we have $\frac{(\{\alpha\}_q;q)_k}{(q;q)_k} = \binom{\alpha + k - 1}{k}_q$. This gives part $(b)$.

    Similarly, for part $(a)$, we can substitute $-\{\alpha\}_qx$ for $x$, and $a = \{\alpha\}_q^{-1}$ to get
    \[ \frac{(-x;q)_\infty}{(-\{\alpha\}_qx;q)_\infty} = \sum_{k=0}^\infty (-1)^k \{\alpha\}_q^k \frac{(\{\alpha\}_q^{-1};q)_k}{(q;q)_k} x^k \]
    By Proposition \ref{prop:alternate_binom_formula}(c), we have $(-1)^k \{\alpha\}_q^k \frac{(\{\alpha\}_q^{-1};q)_k}{(q;q)_k} = q^{\binom{k}{2}} \binom{\alpha}{k}_q$.
\end {proof}

It will be convenient to give names to these series, and so we will call them $B_\alpha$ and $b_\alpha$:

\begin {defn}
    Let $\alpha \in \Bbb{R}$, and define $B_\alpha(q,x)$ and $b_\alpha(q,x)$ as the following elements of $\Bbb{Q}((q))[[x]]$:
    \[ B_\alpha(q,x) := \sum_{k=0}^\infty q^{\binom{k}{2}} \binom{\alpha}{k}_q x^k = \frac{(-x;q)_\infty}{(-\{\alpha\}_qx;q)_\infty} \]
    \[ b_\alpha(q,x) := \sum_{k=0}^\infty \binom{\alpha + k - 1}{k}_q x^k = \frac{(\{\alpha\}_q x; q)_\infty}{(x;q)_\infty} \]
\end {defn}

The following result is a $q$-analogue of the simple statement that $(1+x)^\alpha = (1+x)(1+x)^{\alpha-1}$.

\begin {prop} \label{thm:alpha_plus_one}
    For $\alpha \in \Bbb{R}$ and $n \in \Bbb{Z}$,
    \begin {enumerate}
        \item[$(a)$] $\displaystyle B_{\alpha+1}(q,x) = (1+x) \cdot B_{\alpha}(q, qx) = (1 + \{\alpha\}_qx) \cdot B_\alpha(q,x)$
        \item[$(b)$] $\displaystyle b_{\alpha+1}(q,x) = \frac{1}{1-x} \cdot b_{\alpha}(q, qx) = \frac{1}{1-\{\alpha\}_q x} \cdot b_\alpha(q,x)$
        \item[$(c)$] $\displaystyle B_{\alpha+n}(q,x) = B_n(q,x) \cdot B_{\alpha}(q, q^n x) = B_n(q,\{\alpha\}_q x) \cdot B_\alpha(q,x)$
        \item[$(d)$] $\displaystyle b_{\alpha+n}(q,x) = b_n(q,x) \cdot b_{\alpha}(q, q^n x) = b_n(q,\{\alpha\}_q x) \cdot b_\alpha(q,x)$
    \end {enumerate}
\end {prop}
\begin {proof}
    Parts $(c)$ and $(d)$ obviously follow from $(a)$ and $(b)$ by induction. Parts $(a)$ and $(b)$ follow from the product formula in
    Theorem \ref{thm:product_formula}. Parts $(a)$ and $(b)$ can also be derived by manipulating the power series 
    using the $q$-Pascal identity (Proposition \ref{prop:q_Pascal}).
    \begin {align*} 
        B_\alpha(q,x) &= \sum_{k \geq 0} q^{\binom{k}{2}} \binom{\alpha}{k}_q x^k \\
                 &= \sum_{k \geq 0} q^{\binom{k}{2}}q^k \binom{\alpha-1}{k}_q x^k 
                  + \sum_{k \geq 0} q^{\binom{k}{2}} \binom{\alpha-1}{k-1}_q x^k \\
                 &= \sum_{k \geq 0} q^{\binom{k}{2}}q^k \binom{\alpha-1}{k}_q x^k 
                  + \sum_{k \geq 0} q^{\binom{k}{2}}q^k \binom{\alpha-1}{k}_q x^{k+1} \\
                 &= (1 + x) \sum_{k \geq 0} q^{\binom{k}{2}} \binom{\alpha-1}{k}_q (qx)^k \\
                 &= (1 + x) B_{\alpha-1}(q,qx)
    \end {align*}
    For the second equality in part (a), use Proposition \ref{prop:q_Pascal} (b). The calculations for $b_\alpha$ are similar.
\end {proof}

\bigskip

\begin {ex}
    The power series $b_{1/2}(q,x)$ is a $q$-analogue of the function $\frac{1}{\sqrt{1-x}}$. By definition, we have
    \[ b_{1/2}(q,x) = \sum_{k=0}^\infty \binom{1/2 + k-1}{k}_q x^k \]
    Using the fact that $\left\{ \frac{1}{2} \right\}_q = \frac{1+q^2}{1+q}$, we get fom Proposition \ref{prop:alternate_binom_formula}(e), and some algebraic manipulation, that
    \[ \binom{1/2 + k-1}{k}_q  = \frac{\left( \frac{1+q^2}{1+q}; q \right)_k}{(q;q)_k} = \prod_{i=0}^{k-1} \frac{1+q-q^i-q^{i+2}}{1 + q - q^{i+1} - q^{i+2}} \]
    The first few terms are given by
    \begin {align*}
        b_{1/2}(q,x) &= 1 + \frac{q}{1+q}x + \frac{q(1+q+q^2)}{(1+q)^3} x^2 + \frac{q(1+2q+q^3+q^4)}{(1+q)^4} x^3 \\
        &\phantom{=} + \frac{q(1+4q+7q^2+8q^3+7q^4+5q^5+2q^6+q^7)}{1+6q+16q^2+26q^3+30q^4+26q^5+16q^6+6q^7+q^8} x^4 + \cdots
    \end {align*}
\end {ex}

\bigskip

\section {$q$-Calculus}

\bigskip

\begin {defn}
    For a function $f(x)$, its ``\emph{$q$-derivative}'' is
    \[ D_q f(x) := \frac{f(qx) - f(x)}{qx - x} = \frac{f(qx) - f(x)}{(q-1)x} \]
\end {defn}

The function $B_\alpha(q,x)$ is the $q$-analogue of $(1+x)^\alpha$, which has derivative 
$\frac{d}{dx} B_\alpha(1,x) = \alpha(1+x)^{\alpha-1} = \alpha B_{\alpha-1}(1,x)$. Similarly,
$b_\alpha(1,x) = \frac{1}{(1-x)^\alpha}$ has derivative $\frac{\alpha}{(1-x)^{\alpha+1}} = b_{\alpha+1}(1,x)$.
The folowing are $q$-analogues of these statements:

\bigskip

\begin {prop}
    The $q$-derivatives of $B_\alpha(q,x)$ and $b_\alpha(q,x)$ are given as follows.
    \begin {enumerate}
        \item[$(a)$] $D_q B_\alpha(q,x) = [\alpha]_q B_{\alpha-1}(q,qx)$
        \item[$(b)$] $D_q b_\alpha(q,x) = [\alpha]_q b_{\alpha+1}(q,x)$
    \end {enumerate}
\end {prop}
\begin {proof}
    We only need the identity $D_q x^n = [n]_q x^{n-1}$, and we can differentiate the power series term-by-term:
    \begin {align*}
        D_q B_\alpha(q,x) &= D_q \sum_{n \geq 0} q^{\binom{n}{2}} \binom{\alpha}{n}_q x^n \\
                          &= \sum_{n \geq 0} q^{\binom{n}{2}} \binom{\alpha}{n}_q [n]_q x^{n-1} \\
                          &= \sum_{n \geq 0} q^{\binom{n+1}{2}} \binom{\alpha}{n+1}_q [n+1]_q x^n
    \end {align*}
    At this point, note that $\binom{n+1}{2} = \binom{n}{2} + n$, so we can write $q^{\binom{n+1}{2}} = q^{\binom{n}{2}} \cdot q^n$, and
    it is easy to check that 
    \[ \binom{\alpha}{n+1}_q [n+1]_q = [\alpha]_q \cdot \binom{\alpha-1}{n}_q \]
    Using these two observations, we have that
    \[ D_q B_\alpha(q,x) = [\alpha]_q \sum_{n \geq 0} q^{\binom{n}{2}} \binom{\alpha-1}{n}_q q^n x^n = [\alpha]_q B_{\alpha-1}(q,qx) \]
    The computation for $b_\alpha$ is similar, but does not give an extra factor of $q^n$:
    \begin {align*}
        D_q b_\alpha(q,x) &= D_q \sum_{n \geq 0} \binom{\alpha + n - 1}{n}_q x^n \\
                          &= \sum_{n \geq 0} \binom{\alpha + n - 1}{n}_q [n]_q x^{n-1} \\
                          &= \sum_{n \geq 0} \binom{\alpha + n}{n+1}_q [n+1]_q x^n \\
                          &= [\alpha]_q \sum_{n \geq 0} \binom{\alpha + n}{n}_q x^n \\
                          &= [\alpha]_q b_{\alpha+1}(q,x)
    \end {align*}
\end {proof}

\bigskip

\begin {rmk}
    Combining these formulas with the identities from Theorem \ref{thm:alpha_plus_one}, we obtain the $q$-differential equations
    \[ 
        D_q B_\alpha(q,x) = \frac{[\alpha]_q}{1+x} B_\alpha(q,x) 
        \quad \text{and} \quad 
        D_q b_\alpha(q,x) = \frac{[\alpha]_q}{1-x} b_\alpha(q,qx)
    \]
\end {rmk}

\bigskip

\begin {rmk} 
    Theorem \ref{thm:product_formula} can also be proved using these $q$-derivative expressions.
    Starting with the differential equation above, we can obtain a functional equation satisfied by $B_\alpha(q,x)$:
    \begin {align*}
        \frac{[\alpha]_q}{1+x} B_\alpha(q,x) &= \frac{B_\alpha(q,qx) - B_\alpha(q,x)}{(q-1)x} \\
        [\alpha]_q (q-1)x B_\alpha(q,x) &= (1+x) \Big( B_\alpha(q,qx) - B_\alpha(q,x) \Big) \\
        \left(1 + \left( 1 + [\alpha]_q (q-1) \right) x \right) B_\alpha(q,x) &= (1+x) B_\alpha(q,qx) \\
         B_\alpha(q,x) &= \frac{1+x}{1 + \{\alpha\}_q x} B_\alpha(q,qx)
    \end {align*}
    In the last step, we have used that $\{\alpha\}_q = 1 + (q-1)[\alpha]_q$ (Proposition \ref{prop:q^alpha_properties}(a)).
    Now that we have this functional equation, we can use it again on the right-hand side to replace $B_\alpha(q,qx)$
    with $\frac{1+qx}{1 + \{\alpha + 1\}_q x} B_\alpha(q,q^2x)$. Repeatedly applying the functional equation gives the
    infinite product formula from Theorem \ref{thm:product_formula}.

    The calculation for $b_\alpha(q,x)$ is similar.
\end {rmk}

\section{Some More Identities}

Let us now demonstrate some more identities which hold in our setting.
We first give another $q$-analog of the Pascal identiy for $q$-deformed rational and real numbers.
The identity in the following proposition for $q$-integers appeared in~\cite[Theorem 1]{Shannon}.

\begin{prop}\label{prop:otherPascal}
For $\alpha \in \Bbb{R}$ and $k \in \Bbb{Z}$,
\[\displaystyle \binom{\alpha-1}{k}_q + \binom{\alpha-1}{k-1}_q = \left(\frac{2-q^k-\{\alpha-k\}_q}{1 - \{\alpha\}_q}\right) \binom{\alpha}{k}_q \]
\end{prop}
\begin{proof}
We begin computing
\begin{align*}
\binom{\alpha-1}{k}_q + \binom{\alpha-1}{k-1}_q &= \frac{[\alpha-1]_q [\alpha-2]_q \cdots [\alpha-k]_q}{[k]_q!} + \frac{[\alpha-1]_q [\alpha-2]_q \cdots [\alpha-k+1]_q}{[k-1]_q!}\\
&= \left( \frac{[\alpha-k]_q}{[\alpha]_q} + \frac{[k]_q}{[\alpha]_q}\right) \frac{[\alpha]_q [\alpha-1]_q \cdots [\alpha-k+1]_q}{[k]_q!}\\
&= \left( \frac{(1-\{\alpha-k\}_q) + (1-q^k)}{1-\{\alpha\}_q}\right) \binom{\alpha}{k}_q
\end{align*}
and the result follows noting we have made use of Propostion~\ref{prop:q^alpha_properties} part (b).
\end{proof}

We also have a Chu--Vandermonde identity in the setting of $q$-deformed rational and real numbers.

\begin{cor}[{$q$-Chu--Vandermonde}]\label{cor:Vandermonde}
    For $\alpha \in \Bbb{R}$ and $k,n \in \Bbb{Z}$,
    \[\binom{\alpha + n}{k}_q = \sum_{j = 0}^k q^{j(n-k+j)} \binom{n}{k-j}_q \binom{\alpha}{j}_q \]
\end{cor}
\begin{proof}
    We start with Theorem~\ref{thm:alpha_plus_one} part (c) which states that $B_{\alpha+n}(q,x) = B_n(q,x) \cdot B_{\alpha}(q, q^n x)$.
    As a summation this means
    \[\sum_{k=0}^\infty q^{\binom{k}{2}} \binom{\alpha+n}{k}_q x^k = \left( \sum_{k=0}^n q^{\binom{k}{2}} \binom{n}{k}_q x^k \right) \left( \sum_{k=0}^\infty q^{\binom{k}{2} + nk} \binom{\alpha}{k}_q x^k \right)\]
    from which we will equate the coefficient of $x^k$ on each side. Thus we obtain
    \begin{align*}
        q^{\binom{k}{2}} \binom{\alpha+n}{k}_q &= \sum_{j = 0}^k \left(q^{\binom{k-j}{2}} \binom{n}{k-j}_q \right) \left(q^{\binom{j}{2} + nj} \binom{\alpha}{j}_q\right)\\
        &= \sum_{j=0}^k q^{\binom{k}{2} + j(n-k+j)} \binom{n}{k-j}_q \binom{\alpha}{j}_q
    \end{align*}
    and the corollary immediately follows. 
\end{proof}

Let us now give a lemma that is a consequence of the $q$-Chu--Vandermonde identity.

\begin{lem}
For $\alpha \in \Bbb{R}$ and $\ell,m,n \in \Bbb{Z}$,
\[\sum_{j \geq 0} q^{\ell(j-n+\ell) + j(m-n+j)} \binom{\ell}{n-j}_q \binom{\alpha}{m+j}_q = q^{(m-\ell)(n-\ell)} \binom{\alpha + \ell}{m+n}\]
\label{lem:vand}
\end{lem}

\begin{proof}
We begin by applying Corollary~\ref{cor:Vandermonde} to $\binom{\alpha + \ell}{m+n}_q$ and simplifying to obtain
\begin{align*}
\binom{\alpha + \ell}{m+n} &= \sum_k q^{k(\ell - (m+n) + k)} \binom{\alpha}{k}_q \binom{\ell}{m+n-k}_q\\
&= \sum_j q^{(m+j)(\ell-n+j)} \binom{\alpha}{m+j}_q \binom{\ell}{n-j}_q
\end{align*}
which completes the proof after comparing exponents.
Indeed, it is the case that \[(m-\ell)(n-\ell) + (m+j)(\ell-n+j) = \ell(j-n+\ell) + j(m-n+j)\]
which completes the proof.
\end{proof}

The next proposition in the case of integers is found in Riordan's book~\cite[Section 1.4 Equation (10)]{Riordan}, and a $q$-integer version was given in~\cite[Section 3.9]{thesis}.

\begin{prop}
For $\alpha \in \Bbb{R}$ and $m,n \in \Bbb{Z}$,
\[\binom{\alpha}{m}_q \binom{\alpha}{n}_q = \sum_{\ell \geq 0} q^{(m-\ell)(n-\ell)} \binom{n}{\ell}_q \binom{m}{\ell}_q \binom{\alpha+\ell}{m+n}_q \]
\end{prop}

\begin{proof}
We begin computing making use of the $q$-Chu--Vandermonde identity in Corollary~\ref{cor:Vandermonde} as well as other simplifications.
\begin{align*}
\binom{\alpha}{m}_q \binom{\alpha}{n}_q &= \binom{\alpha}{m}_q \binom{(\alpha-m)+m}{n}_q\\
&= \sum_{j\geq 0} q^{j(m-n+j)} \binom{\alpha}{m}_q \binom{m}{n-j}_q \binom{\alpha-m}{j}_q   \tag{Corollary \ref{cor:Vandermonde}}\\
&= \sum_{j\geq 0} q^{j(m-n+j)} \binom{m+j}{n}_q \binom{n}{j}_q \binom{\alpha}{m+j}_q\\
&= \sum_{j \geq 0} q^{j(m-n+j)} \binom{n}{j}_q \binom{\alpha}{m+j}_q \sum_{\ell \geq 0} q^{\ell(j-n+\ell)} \binom{j}{n-\ell}_q \binom{m}{\ell}_q \tag{Corollary \ref{cor:Vandermonde}}\\
&= \sum_{\ell \geq 0} \sum_{j \geq 0} q^{j(m-n+j) + \ell(j-n+\ell)} \binom{m}{\ell}_q \binom{n}{j}_q \binom{j}{n-\ell}_q  \binom{\alpha}{m+j}_q\\
&= \sum_{\ell \geq 0} \sum_{j \geq 0} q^{j(m-n+j) + \ell(j-n+\ell)} \binom{m}{\ell}_q \binom{n}{\ell}_q \binom{\ell}{n-j}_q \binom{\alpha}{m+j}_q\\
&= \sum_{\ell \geq 0} \binom{m}{\ell}_q \binom{n}{\ell}_q \sum_{j \geq 0} q^{j(m-n+j) + \ell(j-n+\ell)}  \binom{\ell}{n-j}_q  \binom{\alpha}{m+j}_q
\end{align*}
The proof is then completed by using Lemma~\ref{lem:vand}.
Unlabeled steps above are routine algebraic manipulations, and many of these steps are simply rearrangement of the factors in the falling factorials in the $q$-binomials.
\end{proof}

\bigskip

It is well-known that for the classical $q$-binomial coefficients, $\lim_{n \to \infty} \binom{n}{k}_q = \frac{1}{(q;q)_k}$. 
We have the following generalization of this fact:

\begin {prop}
    Let $\alpha_1 \leq \alpha_2 \leq \alpha_3 \leq \cdots$ be any sequence of real numbers for which $\lim_{n \to \infty} \alpha_n = \infty$. Then
    \[ \lim_{n \to \infty} \binom{\alpha_n}{k}_q = \frac{1}{(q;q)_k} \]
\end {prop}
\begin {proof}
If $k=0$, then we immediatey find that both sides are equal to $1$.
So, we may assume that $k \geq 1$.
    By Proposition \ref{prop:alternate_binom_formula}(e), we have
    \[ \binom{\alpha}{k}_q = \frac{(\{\alpha-k+1\}_q; q)_k}{(q;q)_k} = \frac{(q^{1-k}\{\alpha\}_q; q)_k}{(q;q)_k} \]
    The numerator is the product
    \[ (1 - q^{1-k}\{\alpha\}_q) (1 - q^{2-k}\{\alpha\}_q) \cdots (1 - \{\alpha\}_q) \]
    If $N = \lfloor \alpha \rfloor$ is the integer part of $\alpha$, then Proposition \ref{prop:q^alpha_properties} says that
    $\{\alpha\}_q = q^N \{\alpha-N\}_q$. 
Since we are considering $\lim_{n \to \infty} \alpha_n$ and $k$ is constant, we may assume that $N\gg k$.
Therefore the power series for $\{\alpha\}_q$ has no terms of degree less than $N$.
    The product given above is therefore of the form
    \[ 1 + \Big( \text{terms of degree at least $N-k+1$} \Big) \]
    Since the sequence $\alpha_n$ increases without bound, we can eventually find $n$ large enough so that $\alpha_n > N$,
    and so $\lim_{n \to \infty} (q^{1-k}\{\alpha_n\}_q; q)_k = 1$. We therefore have
    \[ \lim_{n \to \infty} \binom{\alpha_n}{k}_q = \lim_{n \to \infty} \frac{(q^{1-k}\{\alpha_n\}_q; q)_k}{(q;q)_k} = \frac{1}{(q;q)_k} \]
\end {proof}

\section {A New $q$-Analogue of the Gamma Fuction}

\bigskip

Note that the $q$-Pochhammer symbol is equal to $(x;q)_k = \frac{(x;q)_\infty}{(q^kx;q)_\infty}$. It is therefore natural (e.g. see the appendix of \cite{gasper_book})
to define a version when $k=\alpha$ is not an integer by the formula $\frac{(x;q)_\infty}{(q^\alpha x;q)_\infty}$. Our basic philosophy
of replacing $q^\alpha$ with the function $\{\alpha\}_q$ leads us to the following definition.

\begin {defn} \label{def:q_pochhammer}
    For $\alpha \in \Bbb{R}$, define the generalized $q$-Pochhammer symbol by
    \[ (x;q)_\alpha := \frac{(x;q)_\infty}{(\{\alpha\}_q x;q)_\infty} = \frac{1}{b_\alpha(q,x)} = B_\alpha(q,-x) \]
\end {defn}

\bigskip

\begin {rmk}
    Using this new notation, Theorem \ref{thm:product_formula} can be re-stated as follows:
    \[ (x;q)_\alpha = \sum_{k=0}^\infty (-1)^k q^{\binom{k}{2}} \binom{\alpha}{k}_q x^k \]
    \[ \frac{1}{(x;q)_\alpha} = \sum_{k=0}^\infty \binom{\alpha+k-1}{k}_q x^k \]
\end {rmk}

\bigskip

For a positive integer $n$, the $q$-factorial is related to the Pochhammer symbol by $[n]_q ! = \frac{(q;q)_n}{(1-q)^n}$. 
Accordingly, the $q$-Gamma function is usually taken to be $\frac{(q;q)_\infty}{(q^x;q)_\infty (1-q)^{x-1}}$.
This leads us to a definition of a new $q$-Gamma function after the following lemma.

\begin{lem}
    The evaluation $B_{\alpha-1}(q,-q)$ gives a well-defined power series if and only if $\alpha \geq 1$.
    \label{lem:well}
\end{lem}

\begin{proof}
    A term in the expansion of $B_{\alpha-1}(q,-q)$ looks like $(-1)^k q^{\binom{k}{2} + k} \binom{\alpha-1}{k}_q$. 
    If $N \leq \alpha - 1 < N+1$  for $N \in \Bbb{Z}$ and $k > \alpha - 1$, then by Lemma~\ref{lem:order} when $N \geq 0$ we have 
    \[\ord\left( (-1)^k q^{\binom{k}{2} + k} \binom{\alpha-1}{k}_q \right) = \binom{k}{2} + k + b - \binom{k-N}{2}\]
    where $b \geq 1$ is as in the lemma.
    In this case the order is a strictly increasing function in $k$, implying each degree has only finitely many terms contributing.
    So, $B_{\alpha-1}(q,-q)$ is well-defined as a formal power series when $\alpha \geq 1$.\\

    However, in the case that $N < 0$
    \[\ord\left( (-1)^k q^{\binom{k}{2} + k} \binom{\alpha-1}{k}_q \right) = \binom{k}{2} + k +Nk - \binom{k}{2} = (N+1)k\]
    which is identically $0$ when $N = -1$.
    For $N < -1$ we find the order approaches $- \infty$ as $k$ increases.
    Thus $B_{\alpha-1}(q,-q)$ is not defined even as a formal Laurent series for $\alpha < 1$.
\end{proof}

\begin {defn}
    Define a function $\Gamma_q \colon \Bbb{R} \setminus \mathbb{Z}_{<0} \to \Bbb{R}((q))$ as follows: For $\alpha \geq 1$ we set
    \[ \Gamma_q(\alpha) := \frac{(q;q)_{\alpha-1}}{(1-q)^{\alpha-1}} = \frac{1}{(1-q)^{\alpha-1} b_{\alpha-1}(q,q)} = \frac{B_{\alpha-1}(q,-q)}{(1-q)^{\alpha-1}} \]
    which is well-defined by Lemma~\ref{lem:well}.
    In this expression, the numerator $(q;q)_{\alpha-1}$ is the one from Definition \ref{def:q_pochhammer}, and
    the denominator $\frac{1}{(1-q)^{\alpha-1}}$ is taken to be the formal power series $\sum_n \binom{\alpha+n-2}{n} q^n$,
    where the binomial coefficients $\binom{\alpha+n-2}{n} \in \Bbb{R}$ are the ordinary ones (not the $q$-analogues).
    For $\alpha < 1$ we define 
    \[\Gamma_q(\alpha) := \frac{\Gamma_q(\alpha+1)}{[\alpha]_q} \]
    which is well-defined since it is obtained through division by a Laurent series.
\label{def:Gamma}
\end {defn}

\bigskip

For integers, the usual $\Gamma$ function satisfies $\Gamma(n+1) = n!$
To avoid confusion with our $q$-Gamma function we have given in Definition~\ref{def:Gamma}, we will denote the classical $q$-Gamma function by $G_q(x)$.   
For the classical $q$-Gamma function, we have $G_q(n+1) = [n]_q!$,
and more generally the well known property that $G_q(x+1) = \frac{1-q^x}{1-q} G_q(x)$. We required an analogous shift property to hold for $\alpha < 1$, but we actually have the corresponding property for our new $q$-Gamma function for all input values.

\begin {prop}
    For $\alpha \in \Bbb{R}$ with $\alpha \geq 1$, we have
    \[ \Gamma_q(\alpha+1) = [\alpha]_q \, \Gamma_q(\alpha) \]
\label{prop:GammaShift}
\end {prop}
\begin {proof}
    Use Proposition \ref{thm:alpha_plus_one} to write $B_\alpha(q,-q)$ as $(1-\{\alpha\}_q) B_{\alpha-1}(q,-q)$, and factor
    $\frac{1}{(1-q)^\alpha}$ as $\frac{1}{1-q} \cdot \frac{1}{(1-q)^{\alpha-1}}$. Then we have that $\Gamma_q(\alpha+1) = \frac{1-\{\alpha\}_q}{1-q} \Gamma_q(\alpha)$,
    and by Proposition \ref{prop:q^alpha_properties}, the factor on the left is equal to $[\alpha]_q$.
\end {proof}

For the classical (integer) $q$-binomial coefficients, they can be written in terms of $q$-factorials and $q$-Pochhammer symbols as
\[ \binom{n}{k}_q = \frac{(q;q)_n}{(q;q)_k (q;q)_{n-k}} = \frac{[n]_q !}{[k]_q! [n-k]_q!} \]
Using the definitions given above, we have similar expressions for the $q$-rational and $q$-real binomial coefficients.

\begin {prop}
    For $\alpha \in \Bbb{R}$, and positive integer $k$, we have
    \[ \binom{\alpha}{k}_q = \frac{(q;q)_\alpha}{(q;q)_k (q;q)_{\alpha-k}} = \frac{\Gamma_q(\alpha+1)}{\Gamma_q(k+1) \Gamma_q(\alpha-k+1)} \]
\end {prop}
\begin {proof}
    Assume $\alpha \geq 1$. The two expressions on the right-hand side are obviously equal by the definition of $\Gamma_q(\alpha)$. To see the first equality, notice that
    \[ \frac{(q;q)_\alpha}{(q;q)_{\alpha-k}} = \frac{(\{\alpha-k+1\}_q;q)_\infty}{(\{\alpha+1\}_q;q)_\infty} = (\{\alpha\}_q,q^{-1})_k \]
    Dividing by $(q;q)_k$ gives $\binom{\alpha}{k}_q$ by Proposition \ref{prop:alternate_binom_formula}(d).\\
    
    In the case that $\alpha < 1$ we need to use our shift property. First notice that
    \[ \left( \frac{[\alpha-k]_q}{[\alpha]_q} \right) \binom{\alpha}{k}_q = \binom{\alpha-1}{k}_q\]
    for $\alpha \neq 0$. Hence, if $\binom{\alpha}{k}_q$ has the form specified by the proposition then we see that $\binom{\alpha-1}{k}_q$ does as well by using Proposition~\ref{prop:GammaShift}.
\end {proof}

\bigskip

Note that $\Gamma_q(\alpha)$ is a Laurent series with real (not integer) coefficients. For example,
\[ 
    \Gamma_q \left( \frac{3}{2} \right) = 1 + \frac{1}{2} q - \frac{5}{8} q^2 - \frac{3}{16} q^3 
                                            + \frac{115}{128} q^4 - \frac{401}{256} q^5 + \frac{2383}{1024} q^6 - \frac{8139}{2048}q^7 + \cdots
\]
For the classical $\Gamma$ function, we have Euler's reflection formula, $\Gamma(x)\Gamma(1-x) = \frac{\pi}{\sin(\pi x)}$.
It turns out that for the $q$-version, this expression is a power series which always has integer coefficients.

\bigskip

\begin {prop} \label{prop:reflection_formula}
	Let $\alpha \in \Bbb{R} \setminus \Bbb{Z}$. Then $\Gamma_q(\alpha) \Gamma_q(1-\alpha) \in \Bbb{Z}((q))$. 
	In other words, it is a Laurent series with integer coefficients.
\end {prop}
\begin {proof}
	There are two cases to consider: when $0 < \alpha < 1$ and when $\alpha > 1$. 

	Consider first the case that $\alpha > 1$, 
	and let $N = \lfloor \alpha \rfloor$ be the integer part of $\alpha$. Then $-N < 1-\alpha < -(N-1)$. Then by definition we have
	\[ \Gamma_q(\alpha) = B_{\alpha-1}(q,-q) b_{\alpha-1}(1,q) \quad \text{and} \quad \Gamma_q(1-\alpha) = \frac{B_{N+1-\alpha}(q,-q)b_{N+1-\alpha}(1,q)}{[1-\alpha]_q[2-\alpha]_q \cdots [N+1-\alpha]_q} \]
	Observe that the only possibility for non-integer coefficients comes from the factors of $b_\alpha(1,q) = \frac{1}{(1-q)^\alpha}$, whose
	coefficients are the scalar non-integer binomial coefficients of the form $\binom{\alpha+k-1}{k}$.
	When we multiply these together in the above expression, we get a factor of
	\[ b_{\alpha-1}(1,q) b_{N+1-\alpha}(1,q) = \frac{1}{(1-q)^{\alpha-1}(1-q)^{N+1-\alpha}} = \frac{1}{(1-q)^N} = \sum_{k=0}^\infty \binom{N+k-1}{k} q^k \]
	In particular, this has integer coefficients. The other factors of $B_\alpha$ and $\frac{1}{[1-\alpha]_q}$, etc, all have integer coefficients.

	\bigskip

	Secondly, there is the case where $0 < \alpha < 1$. In this case both $\alpha$ and $1-\alpha$ are less than $1$, so we have
	\[ \Gamma_q(\alpha) = \frac{\Gamma_q(\alpha+1)}{[\alpha]_q} = \frac{B_\alpha(q,-q)b_\alpha(1,q)}{[\alpha]_q} \quad \text{and} \quad \Gamma_q(1-\alpha) = \frac{\Gamma_q(2-\alpha)}{[1-\alpha]_q} = \frac{B_{1-\alpha}(q,-q)b_{1-\alpha}(1,q)}{[1-\alpha]_q} \]
	As above in the other case, when we multiply these we get a factor of $b_\alpha(1,q)b_{1-\alpha}(1,q) = \frac{1}{1-q}$, which has integer coefficients. The remaining factors all have integer coefficients as well.	
\end {proof}

\bigskip

\begin {ex}
    Above, we gave the first several terms of $\Gamma_q \left( \frac{3}{2} \right)$. We can use this to compute
    \[ 
        \Gamma_q \left( \frac{1}{2} \right) = \frac{\Gamma_q \left( \frac{3}{2} \right)}{\left[\frac{1}{2}\right]_q} 
        = q^{-1} + \frac{3}{2} - \frac{1}{8} q - \frac{13}{16} q^2 + \frac{91}{128}q^3 - \frac{171}{256}q^4 
          + \frac{779}{1024} q^5 - \frac{3373}{2048} q^6 + \cdots 
    \]
    As an example of Proposition \ref{prop:reflection_formula}, we can square this expression to obtain
    \[ 
        \Gamma_q \left( \frac{1}{2} \right)^2 = \Gamma_q \left( \frac{1}{2} \right) \Gamma_q \left( 1 - \frac{1}{2} \right) 
        = q^{-2} + 3q^{-1} + 2 - 2q - q^2 + q^3 - 2q^5 + \cdots 
    \]
    Euler's reflection formula says that $\Gamma \left( \frac{1}{2} \right)^2 = \pi$, so this Laurent series with integer coefficients is, 
    in some sense, a $q$-analogue of the number $\pi$.
\end {ex}

\bigskip

There is another situation where we can obtain Laurent series with integer coefficients, given in the following proposition.

\begin {prop}
    Let $\frac{a}{b} \in \Bbb{Q}$. Then $\Gamma_q \left( \frac{a}{b} \right)^b$ is a Laurent series with integer coefficients.
\end {prop}
\begin {proof}
    The idea is similar to the proof of Proposition \ref{prop:reflection_formula}, but even simpler. The only occurrences of non-integer coefficients
    come from the power series expansion of the factor $\frac{1}{(1-q)^{\frac{a-b}{b}}}$. Raising this to the power $b$ gives $\frac{1}{(1-q)^{a-b}}$,
    which has integer coefficients.
\end {proof}

\bigskip

\begin {ex}
    When $\alpha = \frac{2}{3}$, we have
    \[ \Gamma_q \left( \frac{2}{3} \right) = q^{-1} - \frac{2}{3} - \frac{1}{9} q - \frac{166}{81}q^2 + \frac{803}{243}q^3 - \frac{1553}{729}q^4 + \cdots \]
    Raising this to the third power gives
    \[ \Gamma_q \left( \frac{2}{3} \right)^3 = q^{-3} - 2 q^{-2} + q^{-1} - 6 + 18q - 21q^2 + 27q^3 - 69q^4 + \cdots \]
\end {ex}

\bigskip

\section {Further Questions}

We conclude with some lingering questions which were not addressed in the present paper, 
and which we feel would be interesting diretions for further study.

\begin {ques} 
	What other $q$-analogues have a version which uses $q$-rational or $q$-real numbers?
\end {ques}

In this article we have given some examples where the common $q$-analogue of a number $\alpha$, given by $\frac{1-q^\alpha}{1-q}$, can 
be replaced by Morier-Genoud and Ovsienko's $q$-rationals and $q$-reals, with the latter version retaining many of the desirable properties
of the former. The main observation was that $\{\alpha\}_q = [\alpha+1]_q - [\alpha]_q$ was a good substitute for the expression $q^\alpha$
in these examples. The success of this basic idea and philosophy in the examples given here suggest that maybe this idea can be pushed even further.
It would be interesting to find even more well-known $q$-analogues for which replacing $q^\alpha$ by $\{\alpha\}_q$ gives a new $q$-analogue
(and one which is still interesting and meaningful!).

\bigskip

\begin {ques}
	Do the coefficients of the series representation of $q$-real binomial coefficients have some combinatorial interpretation?
\end {ques}

By definition, $\binom{\alpha}{k}_q$ is a Laurent series with integer coefficients. One might wonder if these integers have some combinatorial meaning.
This question seems difficult however, since the meaning of the integer coefficients of the $q$-reals themselves is still somewhat mysterious.

\bigskip

\begin {ques}
	Is there a more explicit $q$-analogue of the expression $\frac{\pi}{\sin(\pi \alpha)}$ appearing in Euler's reflection formula, which uses 
	the $q$-rational (or $q$-real) number $[\alpha]_q$?
\end {ques}

Euler's reflection formula for the classical $\Gamma$ function says $\Gamma(\alpha)\Gamma(1-\alpha) = \frac{\pi}{\sin(\pi \alpha)}$, and
in Proposition \ref{prop:reflection_formula}, we showed that $\Gamma_q(\alpha)\Gamma_q(1-\alpha)$ is a series with integer coefficients
(despite the fact that $\Gamma_q(\alpha)$ need not have integer coefficients). In what sense are these series a ``good'' $q$-analogue of
the expression $\frac{\pi}{\sin(\pi \alpha)}$? Is there some $q$-deformation of the elements of this expression (the number $\pi$ and
the sine function) which make this explicit?

\bigskip

\section* {Acknowledgments} \label{sec:acknowledgments}
\addcontentsline {toc}{section} {\nameref{sec:acknowledgments}}

We would like to thank Vic Reiner for interesting questions and discussions which led to the idea for this project.
Nicholas Ovenhouse was supported by the Simons Foundation grant 327929.
John Machacek was supported by NSF grant DMS-2039316.

\vfill

\bibliographystyle{alpha}
\bibliography{qrats}

\end {document}